\documentclass[12pt,a4paper]{amsart}
\usepackage{amssymb}
\usepackage{mathrsfs}

\headsep=1cm \numberwithin{equation}{section}
\newtheorem{thm}{Theorem}[section]
\newtheorem{cor}[thm]{Corollary}
\newtheorem{lem}[thm]{Lemma}
\newtheorem{prop}[thm]{Proposition}
\newtheorem{question}[thm]{Question}
\theoremstyle{definition}

\theoremstyle{remark}

\numberwithin{equation}{section}

\newcommand{\nat}{\mathbb N}

\newcommand\Supp{\operatorname{Supp }}
\newcommand\Hom{\operatorname{Hom}}
\newcommand\depth{\operatorname{depth }}
\newcommand\Ht{\operatorname{ht}}
\newcommand\V{\operatorname{V}}
\newcommand{\p}{\frak p }
\newcommand{\q}{\frak q }
\newcommand{\fa}{\frak a }
\newcommand{\fb}{\frak b }
\newcommand{\fm}{\frak m }

\newcommand{\To}{\longrightarrow}

\newcommand\Max{\operatorname{Max}}
\newcommand\Min{\operatorname{Min}}

\newcommand{\fc}{\frak c }

\newcommand\Spec{\operatorname{Spec}}
\newcommand\Ass{\operatorname{Ass}}

\begin{document}
\title[Faltings' Local-global Principle and Annihilator Theorem]{Faltings' Local-global Principle and Annihilator Theorem for the finiteness dimensions}
\author{Mohammad Reza Doustimehr}
\address{Department of Mathematics, University of Tabriz, Tabriz, Iran;
and School of Mathematics, Institute for Research in Fundamental Sciences (IPM), P.O. Box 19395-5746, Tehran, Iran.}
\email{m\b{ }doustimehr@tabrizu.ac.ir}
\email{mdoustimehr@yahoo.com}
\thanks{ 
This research  was  in part supported by a grant from IPM 
\\
E-mail: {\it m\b{ }doustimehr@tabrizu.ac.ir} (Mohammad Reza Doustimehr)}%
\begin{abstract}
Let $R$ be a commutative Noetherian ring, $M$ be a finitely generated $R$-module  and  $n$ be a non-negative integer.
In this article, it is shown that for a positive integer $t$, there is a finitely generated submodule $N_i$ of $H_{\fa}^i(M)$ such that $\dim\Supp H_{\fa}^i(M)/N_i<n$ for all $i<t$ if and only if there is a finitely generated submodule $N_{i,\p}$ of $H_{\fa R_{\p}}^i(M_{\p})$ such that $\dim\Supp H_{\fa R_{\p}}^i(M_{\p})/N_{i,\p}<n$ for all $i<t$. This generalizes Faltings' Local-global Principle
 for the finiteness  of local cohomology modules (Faltings' in  Math. Ann. 255:45-56, 1981). Also, it is shown that whenever $R$ is a homomorphic image of a Gorenstein local ring,  then  the invariants $\inf\{i\in\nat_0\mid\dim\Supp(\fb^tH_{\fa}^i(M))\geqslant n\text{ for all } t\in\nat_0\}$ and  $\inf\{\depth M_{\p}+\Ht(\fa+\p)/\p\mid\p\in\Spec(R)\setminus\V(\fb), \dim R/(\fa+\p)\geqslant n\}$ are equal, for every finitely generated $R$-module $M$ and for all ideals $\frak a, \frak b$ of $R$ with  $\frak b\subseteq \frak a$.
As a consequence, we determine the least integer $i$ where the local cohomology module $H_{\fa}^i(M)$ is not minimax (resp. weakly laskerian).
\end{abstract}
\maketitle
\section{Introduction}
Throughout this paper, let $R$ denote a commutative Noetherian ring (with identity) and $\frak a$ an ideal of $R$. For an $R$-module $M$, the
$i$th local cohomology module of $M$ with support in $V(\frak a)$ is defined as:
$$H^i_{\frak a}(M) = \underset{n\geq1} {\varinjlim}\,\, {\rm Ext}^i_R(R/\frak a^n, M).$$
Local cohomology was first defined and studied by Grothendieck. We refer the reader to \cite{BS} or \cite{Gr1} for more details about
local cohomology.  An important theorem in local cohomology is Faltings' Local-global Principle for the finiteness of
local cohomology modules \cite[Satz 1]{Fa1}, which states that for a positive integer $t$,   the $R_{\frak p}$-module $H^i_{\frak
aR_{\frak p}}(M_{\frak p})$ is finitely generated for all $i< t$ and for all $\p\in\Spec(R)$ if and only if
the $R$-module $H^i_{\fa}(M)$ is finitely generated for all $i< t$.

 In \cite{AN}, Asadollahi and Naghipour introduced the class of  in dimension $<n$ modules. An $R$-module $M$ is said to be {\it in dimension $<n$}, if there is a finitely generated submodule $N$ of $M$ such that $\dim\Supp M/N<n$.

  In section 2, we extend Faltings' Local-global Principle for the class of in dimension $<n$ modules. More precisely, as a first main result of this paper, we prove the following.
 \begin{thm}
 Let $\fa$ be an ideal of $R$,   $M$ be a finitely generated $R$-module, $n$ be a non-negative integer and $t\in\nat$. Then the following statements are equivalent:
\begin{enumerate}
	\item 	
$H_{\fa}^{i}(M)$ is in dimension $<n$	for all $i<t$;
	\item
$H_{\fa R_{\p}}^{i}(M_{p})$ is in dimension $<n$ for all $i<t$ and all $\p\in\Spec(R)$;
	\item
$H_{\fa R_{\fm}}^{i}(M_{\fm})$ is in dimension $<n$ for all $i<t$ and all $\fm\in\Spec(R)$.
\end{enumerate}
 \end{thm}
 As a generalization of the above theorem, we  prove the following.
\begin{thm}
Let $R$ be a homomorphic image of a Gorenstein ring, $\fa$ and $\fb$ be ideals of $R$ such that $\fb\subseteq\fa$, $n$ be a non-negative integer, and $t\in\nat$. Then the following statements are equivalent:
 \begin{enumerate}
	\item
	There exists $s\in\nat$ such that $\mathfrak{b}^{s}H_{\fa}^{i}(M)$  is in dimension $<n$ for all $i<t$;
	\item
	There exists  $s_{\p}\in\nat$ such that $(\fb R_{\p})^{s_{\p}}H_{\fa R_{\p}}^{i}(M_{\p})$  is in dimension $<n$ for all $i<t$ and all $\p\in\Spec(R)$.
\end{enumerate}
\end{thm}
Another important theorem in local cohomology is Faltings' Annihilator Theorem \cite{Fa2} for local cohomology modules, which states that, if $R$ is a homomorphic image of a regular ring or $R$ has a dualizing complex, then the invariants $f_{\fa}^{\fb}(M)$ and $\lambda_{\fa}^{\fb}(M)$ are equal, for every choice of the finitely generated $R$-module $M$ and for every choice of the ideals $\frak a, \frak b$ of $R$ with  $\frak b\subseteq \frak a$, where
 $f_{\fa}^{\fb}(M)=\inf\{i\in\nat_0\mid\fb\not\subseteq {\rm Rad}(0:H_{\fa}^i(M))\}$ (resp.
$\lambda_{\fa}^{\fb}(M)=\inf\{ {\depth} M_{\p}+ {\Ht}(\fa+\p)/{\p}\mid \p\in\Spec(R)\setminus\V(\fb)\}$)
is the $\frak b$-{\it finiteness dimension of $M$ relative to} $\frak a$ (resp.  the  $\frak b$-{\it minimum  $\frak a$-adjusted depth of $M$}), see \cite[Definitions 9.1.5 and 9.2.2]{BS}.

Recently, Khashyarmanesh and Salarian in \cite{KS}, established a  generalization of  the Faltings' Annihilator Theorem over Gorenstein rings.

In \cite{DN}, the  author and Naghipour defined the $n$th $\fb$-{\it finiteness dimension $f_{\fa}^{\fb}(M)_n$ of $M$ relative to} $\fa$ by
\begin{align*}
f_{\fa}^{\fb}(M)_n:=\inf\{i\in\nat_0\mid {\dim\Supp}\,\fb^tH_{\fa}^i(M)\geq n \text{ for all }t\in\nat_0\},
\end{align*}
and proved that if $R$ is a homomorphic image of a Gorenstein ring, then the invariants $f_{\fa}^{\fb}(M)_n$ and
$\lambda_{\fa}^{\fb}(M)_n:=\inf\{\lambda_{\fa R_{\p}}^{\fb R_{\p}}(M_{\p})\mid \p\in\Spec(R), \dim R/\p\geqslant n\}$ are
equal, for every choice of the finitely generated $R$-module $M$ and for every choice of the ideals $\fa, \fb$ of $R$ with  $\fb\subseteq \fa$.

In section 3, we establish a generalization of Faltings' Annihilator Theorem for the finiteness dimensions. More precisely, as a second main result, we prove the following.
\begin{thm}
Assume that $R$ is a homomorphic image of a Gorenstein local ring. Let
$\fa$ and $\fb$ be ideals of $R$ such that $\fb\subseteq\fa$, and
let $M$ be a finitely generated $R$-module. Then, for each $n\in\nat_0$,
\begin{equation*}
f_{\fa}^{\fb}(M)_n=\inf\{ \depth M_{\p}+ {\Ht}(\fa+\p)/{\p}\mid \p\in\Spec(R)\setminus\V(\fb), \dim R/(\fa+\p)\geqslant n\}.
\end{equation*}
\end{thm}
Our proof is rather similar to the proof of \cite[Theorem 9.3.5]{BS}.  As a consequence of this, we prove that
\begin{equation*}
\mu_{\fa}(M)=\inf\{ \depth M_{\p}+ {\Ht}(\fa+\p)/{\p}\mid \p\in\Spec(R)\setminus\V(\fb), \dim R/(\fa+\p)\geqslant1\},
\end{equation*}
and
\begin{equation*}
\omega_{\fa}(M)=\inf\{ \depth M_{\p}+ {\Ht}(\fa+\p)/{\p}\mid \p\in\Spec(R)\setminus\V(\fb), \dim R/(\fa+\p)\geqslant2\},
\end{equation*}
where $\mu_{\fa}(M)$ (resp. $\omega_{\fa}(M)$) is the least integer $i$ such that the local cohomology module $H_{\fa}^i(M)$ is not minimax (resp. weakly laskerian).

Throughout this paper, $R$ will always be a commutative Noetherian
ring with non-zero identity and $\frak a$ will be an ideal of $R$. An $R$-module $L$ is said to be {\it minimax}, if there exists a finitely generated submodule $N$ of $L$, such that $L/N$ is Artinian. The class of minimax modules was introduced by H. Z\"{o}schinger \cite{Zo1} and he has given in \cite{Zo1, Zo2} many equivalent conditions for a module to be minimax.  A module is called {\it  weakly Laskerian}, if each of its homomorphic images has only finitely many associated primes. The class of weakly laskerian modules was introduced by Divaani-Aazar and Mafi \cite{DM}.
We shall use $\Max R$ to denote the set of all maximal ideals of $R$. Also, for any ideal $\frak a$ of $R$, we denote $\{\frak p \in\Spec(R):\, \frak p\supseteq \frak a \}$ by
$V(\frak a)$. For any unexplained notation and terminology we refer the reader to \cite{BS} and \cite{Mat}.
\section{Local-global Principle}
In this section, we extend Faltings' Local-global Principle for the class of {\it in dimension $<n$} modules. The main results of this section are Theorem \ref{thm2.1} and Theorem \ref{thm2.6}.
We begin with the following lemma which is a generalization of \cite[Lemma 3.1]{AM} and we need it in this paper.
\begin{lem}\label{2.0}
Let $\fa$ be an ideal of $R$, $M$ be an $R$-module and $\mathcal{S}$ be a Serre subcategory of $R$-modules and $R$-homomorphisms. Then $\fa M\in\mathcal{S}$ if and only if
 $M/(0:_M \fa)\in\mathcal{S}$.
\end{lem}
\proof
The proof is exactly similar to \cite[Lemma 3.1]{AM}.\qed\\

 For a non-negative integer $n$ and  a subset  $T$ of $\Spec(R)$, we denote $\{\p\in T|\dim R/\p\geqslant n\}$ by $T_{\geqslant n}$.

  We now state and prove one of the main results of this paper.
\begin{thm}\label{thm2.1}
	Let $\fa$ be an ideal of $R$,   $M$ be a finitely generated $R$-module, $n$ be a non-negative integer and $t\in\nat$. Then the following statements are equivalent:
\begin{enumerate}
	\item 	
$H_{\fa}^{i}(M)$ is in dimension $<n$	for all $i<t$;
	\item
$H_{\fa R_{\p}}^{i}(M_{p})$ is in dimension $<n$ for all $i<t$ and all $\p\in\Spec(R)$;
	\item
$H_{\fa R_{\fm}}^{i}(M_{\fm})$ is in dimension $<n$ for all $i<t$ and all $\fm\in\Spec(R)$.
\end{enumerate}
\end{thm}
\proof
$1\Rightarrow 2\Rightarrow 3)$ follow from definition.

$3\Rightarrow 1)$ Let $s:=\inf \{i|H_{\fa}^{i}(M) \text{ is not in dimension }<n\} $. We suppose that $s<t$ and look for a contradiction. Since
$H_{\fa}^{0}(M), H_{\fa}^{1}(M),...,H_{\fa}^{s-1}(M)$ are in dimension $<n$, it follows from \cite[Theorem 2.2]{AT} that $\Hom_R(R/\fa^u,H_{\fa}^{s}(M))$ is in dimension $<n$ for all $u\in\nat$. Hence, by \cite[Lemma 2.6]{MNS},
$$\Ass H_{\fa}^{s}(M)_{\geqslant n} = \Ass\Hom_{R}(R/\fa,H_{\fa}^{s}(M))_{\geqslant n}$$
is a finite set: Let its members be $\p_1,\p_2,...,\p_r$. For each $1\leqslant j\leqslant r$, there is a maximal ideal $\fm_{j}$ of $R$ such that $\dim \Supp R/\p_j=\dim\Supp R_{\fm_j}/\p_jR_{\fm_j}\geqslant n$. By $(3)$, there is a finitely generated submodule $N_j$ of $H_{\fa}^{s}(M)$ such that $\dim\Supp H_{\fa R_{\fm_j}}^{s}(M_{\fm_j})/(N_j)_{\fm_j}<n$. Thus
$$H_{\fa R_{\p_j}}^s(M_{\p_j})/(N_j)_{\p_j}\simeq(H_{\fa R_{\fm_j}}^s(M_{\fm_j})/(N_j)_{\fm_j})_{\p_jR_{\fm_j}}=0.$$
So $H_{\fa R_{\p_j}}^{s}(M_{\p_j})$ is a finitely generated $R_{\p_j}$-module for all $1\leqslant j\leqslant r$. Hence, there is $u_j\in\nat$ such that $(\fa^{u_j}H_{\fa}^{s}(M))_{\p_j}\simeq (\fa R_{\p_j})^{u_{j}}H_{\fa R_{\p_j}}^{s}(M_{\p_j})=0$ for all $1\leqslant j\leqslant r$. Let $u:=\max\{u_1,u_2,...,u_r\}$. Then $(\fa^{u}H_{\fa}^{s}(M))_{\p_j}=0$ for all $1\leqslant j\leqslant r$.  Since $\Ass(\fa^uH_{\fa}^{s}(M))_{\geqslant n}\subseteq \{\p_1,\p_2,...,\p_r\}$,  $\dim\Supp\fa^uH_{\fa}^{s}(M)<n$. Hence, by Lemma \ref{2.0},
$$\dim\Supp(H_{\fa}^{s}(M)/(0:_{H_{\fa}^{s}(M)}\fa^u))<n.$$
 Now, it follows from the exact sequence
$$0\rightarrow (0:_{H_{\fa}^{s}(M)}\fa^u)\rightarrow H_{\fa}^{s}(M)\rightarrow H_{\fa}^{s}(M)/(0:_{H_{\fa}^{s}(M)}\fa^u)\rightarrow 0$$
that $H_{\fa}^{s}(M)$ is in dimension $<n$, which is a contradiction.
\qed\\

Let $n$ be a nonnegative integer. Bahmanpour et al. in \cite{BNS}, defined the notion of the $n$th finiteness dimension $f_{\fa}^n(M)$ of $M$ relative to $\fa$ by $f_{\fa}^n(M)=\inf\{ f_{\fa R_{\p}}(M_{\p})\mid \p\in{\Supp}M/\fa M, \dim R/\p\geqslant n\}$. In the following corollary, we restate the above theorem in terms of the finiteness dimensions.
\begin{cor}	\label{2.2}
	Let $M$ be a finitely generated $R$-module and $n\in\nat_0$. Then
	\begin{align*}
	f_{\fa}^n(M)=&\inf\{f_{\fa R_{\p}}^n(M_{\p})|\p \in\Spec(R)\}\\
	=&\inf\{f_{\fa R_{\fm}}^n(M_{\fm})|\fm \in\Spec(R)\}.
	\end{align*}
\end{cor}
\proof
By \cite[Theorem 2.10]{MNS}, $f^n_{\fa}(M)=\inf \{i|H_{\fa}^{i}(M) \text{ is not in dimension}\ <n\}$. Now, the result follows from Theorem \ref{thm2.1}.
\qed\\

The special case of Theorem \ref{thm2.1} in which $n=0$ is Faltings' Local-global Principle for the finiteness of local cohomology modules.
\begin{thm} (Faltings' Local-global Principle \cite[Satz 1]{Fa1})
	Let $M$ be a finitely generated $R$-module, $\fa$ be an ideal of $R$ and  $t\in\nat$. Then the following statements are equivalent:
	\begin{enumerate}
		\item 	
		$H_{\fa}^{i}(M)$ is finitely generated	for all $i<t$;
		\item
		$H_{\fa R_{\p}}^{i}(M_{\p})$ is finitely generated for all $i<t$ and all $\p\in\Spec(R)$;
		\item
		$H_{\fa R_{\fm}}^{i}(M_{\fm})$ is finitely generated for all $i<t$ and all $\fm\in\Max(R)$.
	\end{enumerate}
\end{thm}
\proof
The result follows from \cite[Remark 2.2(1)]{MNS} and Theorem \ref{thm2.1}.
\qed\\

Also, as a special case of Theorem \ref{thm2.1} in which $n=1$, we conclude  Faltings' Local-global Principle for the minimaxness of local cohomology modules.

\begin{thm}(Faltings' Local-global Principle for the minimaxness)
Let $M$ be a finitely generated $R$-module, $\fa$ be an ideal of $R$ and  $t\in\nat$. Then the following statements are equivalent:
\begin{enumerate}
	\item
	$H_{\fa}^{i}(M)$ is minimax for all $i<t$;
	\item
	$H_{\fa R_{\p}}^{i}(M_{\p})$ is minimax for all $i<t$ and all $\p\in\Spec(R)$;
	\item
	$H_{\fa R_{\fm}}^{i}(M_{\fm})$ is minimax for all $i<t$ and all $\fm\in\Max(R)$.
\end{enumerate}
\end{thm}
\proof
It follows from \cite[Corollary 2.4]{BNS} and \cite[Theorem 2.10]{MNS} that
\begin{equation*}
\inf\{i|H_{\fa}^{i}(M)\text{ is not in dimension} <1\}=\inf\{i|H_{\fa}^{i}(M)\text{ is not minimax}\}.
\end{equation*}
So $H_{\fa}^{i}(M)$ is in dimension $<1$ for all $i<t$ if and only if $H_{\fa}^i(M)$ is minimax for all $i<t$.
The result now follows from Theorem \ref{thm2.1} for $n=1$.  \qed\\


Z. Tang in \cite{T} proved that    $H_{\fa}^i(M)$ is Artinian for all $i<t$ if and only if $H_{\fa R_{\p}}^{i}(M_{\p})$ is Artinian for all $i<t$ and all $\p\in\Spec(R)$. But this version of Local-global Principle is not concluded from Theorem \ref{thm2.1}, for some non-negative integer $n$. So, it is rather natural to ask the following question.
\begin{question}
How can we generalize Theorem \ref{thm2.1} such that we can conclude all versions of  Local-global Principle as a special case of it?
\end{question}
Now, we prepare the ground with the following proposition to state and prove the next main result of this section. We obtained Proposition \ref{prop2.5} and Theorem \ref{thm2.6} a few years  ago and presented them, as a part of a talk, at the {\it 10th seminar on Commutative algebra and Related Topics of School of Mathematics of the Institute for Research in
Fundamental Sciences (IPM)(Tehran, December 18-19, 2013)}.

\begin{prop}\label{prop2.5}
	Let $M$ be a finitely generated $R$-module and $t\in\nat$. Then the following statements are equivalent:
	\begin{enumerate}
		\item
		There exists $s\in\nat_0$ such that $\fa^sH_{\fa}^i(M)$ is in dimension $<n$ for all $i<t$;
		\item
		$H_{\fa}^i(M)$ is in dimension $<n$.
	\end{enumerate}
\end{prop}
\proof
We use induction on $t$. When $t=1$ , there is  nothing to prove. So suppose that $t>1$ and  the result has been proved for smaller values of $t$.
By the inductive hypothesis $H_{\fa}^{0}(M),H_{\fa}^{1}(M),...,H_{\fa}^{t-2}(M)$ are in dimension $<n$. We must show that $H_{\fa}^{t-1}(M)$ is in dimension $<n$.
 Since $\fa^{s}H_{\fa}^{t-1}(M)$ is in dimension $<n$, it follows from Lemma \ref{2.0} that $H_{\fa}^{t-1}(M)/(0:_{H_{\fa}^{t-1}(M)}\fa^{s})$ is in dimension $<n$. As $H_{\fa}^{0}(M),H_{\fa}^{1}(M),...,H_{\fa}^{t-2}(M)$ are in dimension $<n$, so  $(0:_{H_{\fa}^{t-1}(M)}\fa^{s})$ is in dimension $<n$. Thus, it follows from the exact sequence
$$0\rightarrow (0:_{H_{\fa}^{t-1}(M)}\fa^{s})\rightarrow H_{\fa}^{t-1}(M)\rightarrow H_{\fa}^{t-1}(M)/(0:_{H_{\fa}^{t-1}(M)}\fa^s)\rightarrow 0$$
that the  $R$-module $H_{\fa}^{t-1}(M)$ is in dimension $<n$. This completes the inductive step.
\qed \\

By Theorem \ref{thm2.1} and Proposition \ref{prop2.5}, it is easy to see that there exists $s\in\nat$ such that  $\fa^{s}H_{\fa}^0(M),\fa^{s}H_{\fa}^1(M),...,\fa^sH_{\fa}^{t-1}(M)$ are in dimension $<n$   if and only if there exists $s_{\p}\in\nat$ such that $(\fa R_{\p})^{s_{\p}}H_{\fa R_{\p}}^0(M_{\p}),(\fa R_{\p})^{s_{\p}}H_{\fa R_{\p}}^1(M_{\p}),...,(\fa R_{\p})^{s_{\p}}H_{\fa R_{\p}}^{t-1}(M_{\p})$ are in dimension $<n$ for all $\p\in\Spec(R)$. Now let $\fb$ be a second ideal of $R$. It is rather natural to ask whether the following statements are equivalent:
\begin{enumerate}
	\item
	There exists  $s\in\nat$ such that $\mathfrak{b}^{s}H_{\fa}^{i}(M)$  is in dimension $<n$ for all $i<t$;
	\item
	There exists  $s_{\p}\in\nat$ such that $(\fb R_{\p})^{s_{\p}}H_{\fa R_{\p}}^{i}(M_{\p})$  is in dimension $<n$ for all $i<t$ and all $\p\in\Spec(R)$.
\end{enumerate}

In the following theorem we give an affirmative answer to this question when $R$ is a homomorphic image of a Gorenstein ring and $\fb\subseteq\fa$.
\begin{thm}\label{thm2.6}
Assume that $R$ is a homomorphic image of a Gorenstein ring. Let $\fa$ and $\fb$ be ideals of $R$ such that $\fb\subseteq\fa$, $n$ be a non-negative integer, and $t\in\nat$. Then the following statements are
equivalent:
\begin{enumerate}
\item
There exists  $s\in\nat_0$ such that  ${\fb}^sH_{\fa}^i(M)$ is in dimension $<n$ for all $i<t$;
\item
There exists  $s_{\p}\in\nat_0$ such that $(\fb R_{\p})^{s_{\p}}H_{\fa R_{\p}}^{i}(M_{\p})$  is in dimension $<n$ for all $i<t$ and all $\p\in\Spec(R)$;
\item
$f_{{\fa}R_\p}^{{\fb}R_\p}(M_{\p})\geqslant t$ for each $\p\in{\Spec}(R)$ with ${\dim} R/\p\geqslant n$.
\end{enumerate}
\end{thm}
\proof $1\Rightarrow2$)
For each $i<t$, let $N_i$ be a finitely generated submodule of $\fb^sH_{\fa}^i(M)$ such that ${\dim\Supp}\fb^sH_{\fa}^i(M)/N_i<n$. Let $\p\in\Spec(R)$.
As
$${\dim\Supp}(\fb R_{\p})^sH_{\fa R_{\p}}^i(M_{\p})/(N_i)_{\p}\leqslant{\dim\Supp}\fb^sH_{\fa}^i(M)/N_i<n,$$
it follows that  $(\fb R_{\p})^sH_{\fa R_{\p}}^i(M_{\p})$ is in dimension $<n$.

$2\Rightarrow3$) Let $\p\in\Spec(R)$ with $\dim R/\p\geqslant n$. There is a maximal ideal $\fm$ of $R$ such that $\p\subseteq\fm$ and $\dim R_{\fm}/\p R_{\fm}\geqslant n$. By (2), there exists $s\in\nat_0$ and a finitely generated submodule $N_i$ of $\fb^sH_{\fa}^i(M)$ such that
${\dim\Supp}(\fb R_{\fm})^sH_{\fa R_{\fm}}^i(M_{\fm})/(N_i)_{\fm}={\dim\Supp}(\fb^sH_{\fa}^i(M)/N_i)_{\fm}<n$, for all $i<t$. So,
$(\fb^sH_{\fa}^i(M))_{\p}\simeq((\fb^sH_{\fa}^i(M))_{\fm})_{\p R_{\fm}}$ is finitely generated as $R_{\p}$-module. Since
$(\fb^sH_{\fa}^i(M))_{\p}$ is $\fa R_{\p}$-torsion,  there exists $s'\in\nat$ such that $(\fb^{s+s'}H_{\fa}^i(M))_{\p}\subseteq(\fb^s\fa^{s'}H_{\fa}^i(M))_{\p}=0$. Therefore,
$f_{{\fa}R_\p}^{{\fb}R_\p}(M_{\p})\geqslant t$.

$3\Rightarrow1$) Now, let $R$ be a homomorphic image of a Gorenstein ring. Then, by \cite[Theorem 2.14]{DN} and (3),
\begin{align*}
f_{\fa}^{\fb}(M)_n & =\inf\{\lambda_{\fa
R_{\p}}^{\fb R_{\p}}(M_{\p})\mid\p\in{\Spec}(R) , {\dim} R/\p\geqslant n\} \\
& = \inf\{f_{\fa R_{\p}}^{\fb R_{\p}}(M_{\p})\mid\p\in\Spec(R) , {\dim} R/\p \geqslant n\} \\
& \geqslant t.
\end{align*}
 So, there exists $s\in\nat_0$ such that for each $i<t$, ${\dim\Supp}\fb^sH_{\fa}^i(M)<n$. Hence, $\fb^sH_{\fa}^i(M)$ is in dimension $<n$ for all $i<t$. \qed\\

Theorem \ref{thm2.6} generalizes some of the results of \cite{BRS}.
\section{Faltings' annihilator theorem for Finitness dimensions}
In this section, we prove a generalization of Faltings' Annihilator Theorem \cite{Fa2} over a ring that is homomorphic image of a Gorenstain local ring.
We first prove in the following proposition that if $R$ is a Gorenstein local ring, then $f_{\fa}^{\fb}(M)_n=\gamma_{\fa}^{\fb}(M)_n$,   where $\gamma_{\fa}^{\fb}(M)_n:=\inf\{ \depth M_{\p}+ {\Ht}(\fa+\p)/{\p}\mid \p\in{\Spec}R\setminus\V(\fb), \dim R/(\fa+\p)\geqslant n\}$.

\begin{prop}\label{thm1}
Assume that $R$ is a Gorenstein local ring. Let $\fa$ and $\fb$ be ideals of $R$ such that $\fb\subseteq\fa$, and let $M$ be a finitely generated $R$-module. Then, for every non-negative integer $n$, $$f_{\fa}^{\fb}(M)_n=\gamma_{\fa}^{\fb}(M)_n.$$
\end{prop}

\proof
If $n=0$, it follows from \cite[Theorem 9.5.1]{BS} that
\begin{align*}
f_{\fa}^{\fb}(M)_0&=f_{\fa}^{\fb}(M)\\
&=\inf\{\depth M_{\p}+\Ht(\fa+\p)/\p|\p\in\Spec(R) \backslash\V(\fb)\}\\
&=\inf\{\depth M_{\p}+\Ht(\fa+\p)/\p|\p\in\Spec(R) \backslash\V(\fb)\,,\,\dim R/(\fa+\p)\geqslant0\}\\
&=\gamma_{\fa}^{\fb}(M)_0.
\end{align*}
Suppose that $n>0$. We first prove that $f_{\fa}^{\fb}(M)_n\leqslant\gamma_{\fa}^{\fb}(M)_n$.
Set $\gamma:=\gamma_{\fa}^{\fb}(M)_n$. If $\gamma=\infty$, then there is nothing to prove; we therefore suppose that $\gamma$ is finite. Thus there exists  $\p\in\Spec(R)\setminus\V(\fb)$ such that $\dim R/(\fa+\p)\geqslant n$ and $\depth M_{\p}+\Ht(\fa+\p)/\p=\gamma$. Also, choose $a\in\fb\setminus\p$.
Let $\q,\q'$ be minimal primes of $\fa+\p$ such that $\dim R/\q=\dim R/(\fa+\p)$ and $\Ht\q'/\p=\Ht(\fa+\p)/\p$. Then $\p$ belongs to the set
$$\Sigma:=\{\delta'\in\Spec(R):\p\subseteq\delta'\subseteq\q\text{ and } a\not\in\delta'\}.$$
Let $\delta$ be a maximal member of $\Sigma$. Now $a\in\fb\subseteq\fa\subseteq\q$, and so $\delta\varsubsetneq\q$. By \cite[Lemma 9.3.4]{BS}, $\Ht \q/\delta=1$. Note that the fact that $a\not\in\delta$ ensures that $\delta\in\Spec(R)\setminus\V(\fb)$. Since $\dim R/\q'\leqslant\dim R/\q=\dim R/(\fa+\delta)$, we can deduce from the definition of $\gamma$, \cite[Lemma 9.3.2]{BS} and \cite[Theorem 31.4]{Mat} that
\begin{align*}
\gamma &\leqslant \depth M_{\delta}+\Ht(\fa+\delta)/\delta\\
&\leqslant \depth M_{\delta}+\dim R/\delta-\dim R/(\fa+\delta)\\
&\leqslant\depth M_{\p}+\Ht \delta/\p+\dim R/\delta-\dim R/(\fa+\delta)\\
&=\depth M_{\p}+\dim R/\p-\dim R/(\fa+\delta)\\
&\leqslant\depth M_{\p}+\dim R/\p-\dim R/\q'\\
&=\depth M_{\p}+\Ht \q'/\p\\
&=\depth M_{\p}+\Ht(\fa+\p)/\p\\
&=\gamma.
\end{align*}
Therefore, $\depth M_{\delta}+\Ht(\fa+\delta)/\delta=\gamma$. It is immediate from $\q\in\Min(\fa+\delta)$ and $\Ht \q/\delta=1$ that $\Ht \q/\delta=\Ht(\fa+\delta)/\delta=1$. Thus, we can replace $\p$ by $\delta$ and so make the additional assumption that $\Ht \q/\p=\Ht(\fa+\p)/\p=1$. Now it is easy to see that
\begin{equation*}
1=\Ht (\fa+\p)/\p\leqslant\Ht (\fa R_{\q}+\p R_{\q})/p R_{\q}\leqslant\Ht \q R_{\q}/\p R_{\q}=1.
\end{equation*}
Hence,
\begin{equation*}
\depth (M_{\q})_{\p R_{\q}}+\Ht(\fa R_{\q}+\p R_{\q})/\p R_{\q}=\depth M_{\p}+\Ht (\fa+\p)/\p=\gamma.
\end{equation*}
So, $\lambda_{\fa R_{\q}}^{\fb R_{\q}}(M_{\q})\leqslant\gamma$.
As $\dim R/\q\geqslant n$, it therefore follows from \cite[Theorem 2.10]{DN} and the definition of $\lambda_{\fa}^{\fb}(M)_n$ that
  $$f_{\fa}^{\fb}(M)_n=\lambda_{\fa}^{\fb}(M)_n\leqslant \lambda_{\fa R_{\q}}^{\fb R_{\q}}(M_{\q})\leqslant\gamma.$$

Now, we prove the converse inequality. By \cite[Theorem 2.10]{DN}, $f_{\fa}^{\fb}(M)_n=\lambda_{\fa}^{\fb}(M)_n$. So it is enough to show that $\lambda_{\fa}^{\fb}(M)_n\geqslant\gamma_{\fa}^{\fb}(M)_n$.  If $\lambda_{\fa}^{\fb}(M)_n=\infty$, there is nothing to prove. We therefore suppose that $\lambda_{\fa}^{\fb}(M)_n=t<\infty$. Hence, there exists $\q\in\Spec (R)$ such that $\dim R/\q\geqslant n$ and $\lambda_{\fa R_{\q}}^{\fb R_{\q}}(M_{\q})=t$. Thus there exists $\p\varsubsetneq\q$ such that
$$\depth M_{\p}+\Ht(\fa+\p)/\p\leqslant\depth (M_{\q})_{\p R_{\q}}+\Ht(\p R_{\q}+\fa R_{\q})/\p R_{q}=t.$$
Since $\dim R/(\fa+\p)\geqslant n$,  $\gamma_{\fa}^{\fb}(M)_n\leqslant t$, as required. \qed\\

The next lemma generalizes \cite[Lemma 9.2.6]{BS}.
\begin{lem}\label{lem13}
Let $\fa$ and $\fb$ be ideals of $R$ such that $\fb\subseteq\fa$. Let $M$ be a finitely generated $R$-module, and let $\fc$ be an
ideal of $R$ such that $\fc\subseteq (0:M)$. Then
$$\gamma_{\fa}^{\fb}(M)_n=\gamma_{(\fa+\fc)/\fc}^{(\fb+\fc)/\fc}(M)_n.$$
\end{lem}

\proof The proof is similar to the proof of \cite[lemma 9.2.6]{BS}.  \qed \\

We are now ready to state and prove the second main result of this paper which is a generalization of the Faltings' Theorem for the annihilation of local cohomology modules  over a ring $R$ that  is a homomorphic image of a  Gorenstein local ring.

\begin{thm}\label{thm14}
Assume that $R$ is a homomorphic image of a Gorenstein local ring. Let
$\fa$ and $\fb$ be ideals of $R$ such that $\fb\subseteq\fa$, and
let $M$ be a finitely generated $R$-module. Then, for each $n\in\nat_0$,
$$f_{\fa}^{\fb}(M)_n=\gamma_{\fa}^{\fb}(M)_n.$$
\end{thm}

\proof By assumption, there is a Gorenstein local ring $R'$ and a surjective homomorphism of Noetherian rings $f:R'\To R$. Let $\fa'$ and $\fb'$ be ideals of $R'$ such that $\fa=\fa' R$ and $\fb=\fb'R$. Then by \cite[Lemma 2.12]{DN}, Lemma \ref{lem13}, and Proposition \ref{thm1},
$$f_{\fa}^{\fb}(M)_n=f_{\fa' R}^{\fb'R}(M)_n=f_{\fa'}^{\fb'}(M)_n=\gamma_{\fa'}^{\fb'}(M)_n=\gamma_{\fa}^{\fb}(M)_n,$$
as required. \qed \\

 For determining the least integer $i$ where the local cohomology module $H_{\fa}^i(M)$ is not minimax (resp. weakly laskerian), we need the following proposition.
\begin{prop}
Let $\fa$ be an ideal of $R$, $M$ be a finitely generated $R$-module and $n$ be a non-negative integer. Then
$$f_{\fa}^{\fa}(M)_n=f_{\fa}^n(M).$$
\end{prop}

\proof
For each $i<f_{\fa}^{\fa}(M)_n$, there exists  $t\in\nat$ such that $\dim\Supp \fa^tH_{\fa}^i(M)<n$. So, for all $\p\in\Spec(R)$ with $\dim R/\p\geqslant n$, $(\fa^tH_{\fa}^i(M))_{\p}=0$. Thus, by \cite[Proposition 9.1.2]{BS},
$$f_{\fa}^{\fa}(M)_n\leqslant\inf\{f_{\fa R_{\p}}(M_{\p})\mid \p\in\Spec(R), \dim R/\p\geqslant n\}=f_{\fa}^n(M).$$

Now, we prove the converse inequality. Let $i<f_{\fa}^n(M)$. By \cite[Theorem 2.10]{MNS}, $H_{\fa}^i(M)$ is in dimension $<n$. Hence, there exists a finitely generated submodule $N_i$ of $H_{\fa}^i(M)$  such that $\dim\Supp H_{\fa}^i(M)/N_i<n$. Since $N_i$ is finitely generated, there is  $t\in\nat$ such that $\fa^tN_i=0$, and so $N_i\subseteq (0:_{H_{\fa}^i(M)}\fa^t)$. It therefore follows from the exact sequence
$$H_{\fa}^i(M)/N_i\To H_{\fa}^i(M)/(0:_{H_{\fa}^i(M)}\fa^t)\To0$$
that $\dim\Supp H_{\fa}^i(M)/(0:_{H_{\fa}^i(M)}\fa^t)<n$. Hence, for all prime ideal $\p$ with $\dim R/\p\geqslant n$,
$$(\fa^tH_{\fa}^i(M))_{\p}=(\fa^t(0:_{H_{\fa}^i(M)}\fa^t))_{\p}=0.$$
Thus $\dim\Supp\fa^tH_{\fa}^i(M)<n$. Therefore, $f_{\fa}^n(M)\leqslant f_{\fa}^{\fa}(M)_n$.
\qed\\

\begin{cor}
Assume that $R$ is a homomorphic image of a Gorenstein local ring. Let
$\fa$  be an ideal of $R$, and let $M$ be a finitely generated $R$-module. Then, $$\mu_{\fa}(M)=\inf\{ \depth M_{\p}+ {\Ht}(\fa+\p)/{\p}\mid \p\in\Spec(R)\setminus\V(\fa), \dim R/(\fa+\p)\geqslant1\},$$
where $\mu_{\fa}(M)=\inf\{ i\mid H_{\fa}^i(M) \text{ is not minimax }\}$.

\end{cor}
\proof Since $f_{\fa}^{\fa}(M)_1=f_{\fa}^1(M)$,  the result follows from Theorem \ref{thm14} and \cite[Corollary 2.4]{BNS}. \qed\\
\begin{cor}
Assume that a semi-local ring $R$ is a homomorphic image of a Gorenstein local ring. Let
$\fa$  be an ideal of $R$, and
let $M$ be a finitely generated $R$-module. Then, $$\omega_{\fa}(M)=\inf\{ \depth M_{\p}+ {\Ht}(\fa+\p)/{\p}\mid \p\in\Spec(R)\setminus\V(\fa), \dim R/(\fa+\p)\geqslant2\},$$
where $\omega_{\fa}(M)=\inf\{ i\mid H_{\fa}^i(M) \text{ is not weakly laskerian }\}$.
\end{cor}
\proof Since $f_{\fa}^{\fa}(M)_2=f_{\fa}^2(M)$,  the result follows from Theorem \ref{thm14} and \cite[Proposition 3.7]{BNS}. \qed\\


\end{document}